\documentclass[11pt]{article}

\usepackage{amssymb,latexsym,amsmath,amsfonts}
\usepackage{subfigure}
\usepackage{graphicx, epsf}
\usepackage[active]{srcltx}

\numberwithin{equation}{section}
\newcommand{\bee}{\begin{equation}}
\newcommand{\ene}{\end{equation}}

\newcommand{\la}{\lambda}

\begin{document}

\title{Convergent Analytic Solutions for Homoclinic Orbits in Reversible and Non-reversible Systems}
\author{S.R. Choudhury\footnote{Department of Mathematics, University of Central Florida, USA, choudhur@cs.ucf.edu} $\;$
G. Gambino\footnote{Department of Mathematics, University of Palermo, Italy, gaetana@math.unipa.it} $\;$}

\maketitle


\begin{abstract}

In this paper, convergent, multi-infinite, series solutions are derived for the homoclinic orbits
of a canonical fourth-order ODE system, in both reversible and non-reversible cases. This ODE
includes traveling-wave reductions of many important nonlinear PDEs or PDE systems, for which
these analytical solutions would correspond to regular or localized pulses of the PDE.
As such, the homoclinic solutions derived here are clearly topical, and they
are shown to match closely to earlier results obtained by homoclinic numerical shooting. In addition, the results for the non-reversible case
go beyond those that have been typically
considered in analyses conducted within bifurcation-theoretic settings.

We also comment on generalizing the treatment here to parameter regimes where solutions homoclinic to exponentially small periodic orbits
are known to exist, as well as another possible extension placing the solutions derived here
within the framework of a comprehensive categorization of ALL possible traveling-wave solutions, both smooth and non-smooth, for our governing ODE.
\end{abstract}

%

\section{Introduction}

In this paper, we initiate a novel analytic approach to the calculation of homoclinic orbits in both reversible
and non-reversible ODE systems.

Homoclinic orbits of dynamical systems have been widely treated in recent years by a variety of approaches. For
instance, an early review integrating bifurcation theoretic and numerical approaches was given in \cite{C98}.
Homoclinic orbits are important in applications for a variety of reasons. In
the context of ODE systems, they are often \textit{anchors} for the local dynamics
in their vicinity. Under certain conditions, their existence may indicate the existence of chaos in their neighborhood
(see \cite{R3} - \cite{R4} for instance).

In a totally different setting, if the governing dynamical system is the traveling-wave ODE for a partial differential equation
or equations, its homoclinic orbits correspond to the solitary wave or \textit{pulse} solutions of the PDE(s), which have
many important uses and applications in nonlinear wave propagation theory, nonlinear optics, and in various other settings
(see \cite{R5} - \cite{R6} for instance).

Given their importance, we will develop an analytic solution which will
cover, or be pertinent to, diverse systems of interest in the relevant regime of parameters.
Towards this end, we start with an important, canonical ODE system which
may, among other applications, result from traveling-wave reductions of a variety of physically-important
nonlinear PDEs (henceforth written NLPDEs). As such, this ODE has served as
an important workhorse in treatments of homoclinic orbits in reversible systems \cite{C98}. We shall
also generalize our treatment and include non-reversible terms (nonlinearities) which
have been considered only rarely, if at all, in earlier work.

This generalized canonical ODE is introduced in Section 2, and followed in Section 3 by a review
of some known results for its homoclinic solutions
in four complementary regions spanning all of the parameter space.

Of these regions, two contain countably infinite numbers, or a
so-called \textit{plethora}, of homoclinic orbits.
In one of these two regions, the orbits are homoclinic to the zero solution.
We derive explicit analytic solutions for the homoclinic orbits
in that region in Sections 4 and 5 for the reversible and non-reversible cases
respectively. If the ODE is a traveling-wave reduction of some PDE, these solutions would correspond to genuine or localized pulses of the PDE.

In the second of the above regions containing an infinity of such orbits, the homoclinic solutions are
known to be homoclinic or asymptotic to
exponentially small periodic orbits. These correspond to \textit{delocalized} solitary traveling wave solutions (with exponentially small tails) of
any corresponding PDE(s). We do not treat the solutions in this second region in this paper, but defer them to a future
report.

The remainder of the paper considers our analytic homoclinic solutions, which have the form of multi-infinite series, in greater detail. So as not to affect the readability of the paper, we defer a convergence proof for our series solutions to Section 6.
Before that, the solutions are numerically summed and plotted in Sections 4 and 5 -- fortunately (and
this is subsequently proven using an inductive argument as part of the convergence proof of Section 6), the later coefficients drop off sharply in magnitude (which is by no means guaranteed \textit{a priori} for such convergent series). The plots in Section 4 are in very close agreement with earlier numerical results for the reversible case
obtained by homoclinic shooting. As for the corresponding homoclinic orbits derived for the non-reversible case in Section 5, to our knowledge, the results
there go beyond most which exist in the literature.

For completeness, we should mention at this point that the only OTHER existing results (see \cite{R5}-\cite{R6} for instance) for homoclinic orbits of our ODE, besides the numerical solutions mentioned above, are from the use of the bifurcation-tracking software AUTO. Although
it would be straightforward to compare our series solutions to those as well, we do not consider them as they have
been numerically benchmarked earlier against those obtained by homoclinic shooting, which we do compare our results against.

Section 6 rounds out the paper with a proof of the convergence of our analytical multi-infinite series solutions for the homoclinic orbits, while
Section 7 summarizes the results, and points out directions for further investigations. In particular,
we make some comments on (non-trivially) extending our treatment to the second region
mentioned above where the solutions are homoclinic to exponentially small periodic orbits. We also remark briefly on a forthcoming paper where the solutions and results here will be juxtaposed with, and placed in the framework of, a comprehensive categorization of ALL possible traveling-wave solutions, both smooth and non-smooth, for our governing ODE.

\section{A canonical fourth-order equation and related physical PDEs}
\setcounter{figure}{0}
\setcounter{equation}{0}

Let us consider the equation:

\begin{equation}\label{eq_genC}
u^{''''}-b u^{''}+a u=f(u, u^{'}, u^{''}, u^{'''}),
\end{equation}

\noindent where $a$ and $b$ are real parameters, $f$ is a nonlinear function, and the $'$ indicates the derivatives with respect to $z$.

As one possible motivation, recall that substituting traveling wave solutions of the following form:

\begin{equation}\label{tr_wave}
u(x,t)=u(x-at)\equiv u(z)
\end{equation}

\noindent into the fifth-order KdV5 or FKdV equation:

\begin{equation}\label{KdV5}
u_t+6uu_x+u_{xxx}+u_{xxxxx}+10uu_{xxx}+20u_xu_{xx}+30u^2u_x=0
\end{equation}

\noindent as done, for instance, in \cite{JY01}, one obtains the following equation:

\begin{equation}\label{KdV5_tr}
u^{''''}+u^{''}-au+3 u^2+10uu^{''}+5u^{'2}+10u^3=0
\end{equation}

\noindent which belongs to the class \eqref{eq_genC}. Numerous other
NLPDEs of importance in various application areas and having
traveling wave-reduced equations of the form \eqref{eq_genC} are discussed in \cite{C98} and the
references there.

\section{Possible solitary wave domains via bifurcation analysis}

\begin{figure}
\begin{center}

\includegraphics[scale=.4]{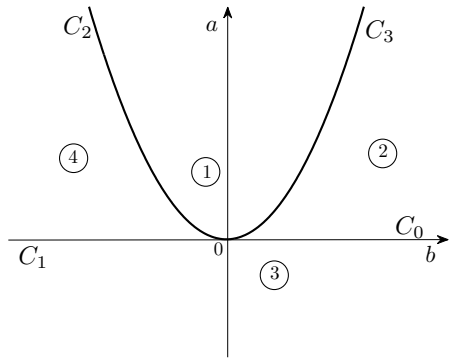}
\end{center}
\caption{\label{parabola} Linearization at the origin of \eqref{eq_genC}.}
\end{figure}

In order to examine the structure of homoclinic solutions to \eqref{eq_genC} the first step is to review its linearization \cite{C98}, which is depicted in Fig.\ref{parabola}.

 In Fig.\ref{parabola} there are four distinct regions, bounded by the codimension one curves $C_i,\, i=0,\dots, 3$, and corresponding to qualitatively different linear dynamics in each part.

Equation \eqref{eq_genC}  is invariant under the transformation $z\to -z$ and is thus a reversible system. In this section, we shall use the theory of reversible systems to characterize the homoclinic orbits to the fixed point of \eqref{eq_genC} in various regions of the $(b,a)$ plane. These would correspond to pulses or solitary wave solutions of any NLPDE (such as those discussed above)
for which \eqref{eq_genC} is the traveling-wave ODE, .

The linearized system corresponding to \eqref{eq_genC}
\bee
u^{''''}-b u^{''}+a u=0\label{lin_eq}
\ene
has a fixed point
\bee
u = u^{'} = u^{''} = u^{'''}=0\label{fixed_pt}.
\ene
Solutions $\phi =ke^{\la z}$ satisfy the characteristic equation $\la^4 -b\la^2 +a=0$ from which one may deduce that the structure of the eigenvalues is distinct in four different regions of the $(b,a)$ plane. These are the regions shown in Figure \ref{parabola} and labeled regions 1-4.

In delineating the structure of the homoclinic orbit to fixed point \eqref{fixed_pt} in various parts of the $(b,a)$ space, we shall first consider the bounding curves $C_0-C_3$ and their neighborhoods. Following  this, we shall discuss the possible occurrence and multiplicities of homoclinic orbits to \eqref{fixed_pt} in each of regions 1 through 4:
\begin{itemize}
\item[a.] {\bf Near $C_0$:} This curve, on which the eigenvalues have the structure $\la_{1-4} =0,0,\pm \la$, and its vicinity have been considered in the context of reversible systems in \cite{ik92}-\cite{kk88}. In this region, a standard analysis yields the normal form on the center manifold
\begin{align*}
\dot x_1 &= x_2\\
\dot x_2 &= \text{sign} (\mu )x_1-\frac32 x_1^2
\end{align*}
where $\mu$ is an unfolding parameter \cite{ik92}. For $\mu>0$, this yields a unique, symmetric homoclinic solution
$$
x_1(t)={\rm sech}^2(t/2)
$$
in the vicinity of $C_0$. One may also show persistence of this homoclinic solution in the original system (2.1) for $\mu >0$ \cite{ik92}.

\item[b.] {\bf Near $C_1$:} Near $C_1$, which corresponds to the eigenvalue structure $\la_{1-4}=0,0,\pm i \omega$, we will show by analysis of a four-dimensional normal form \cite{ik92} in Section 4 that on the side of $C_1$ corresponding to region 3 in Figure 1 there is a ${\rm sech}^2$ homoclinic orbit.

However, in region 3, where the eigenvalue structure is that of a saddle-center $\la_{1-4}=\pm \la,  i\omega$, the fixed point \eqref{fixed_pt} is non-hyperbolic. It may be proved \cite{el96}-\cite{el97} that there are classes of homoclinic orbits in region 3 which are homoclinic to periodic orbits as  $z\to \pm \infty$. These are so-called delocalized solitary waves \cite{jpb91}. Depending on the form of the nonlinear term (the RHS of (2.1) in our case), these periodic orbits often have exponentially small amplitudes. Also, on isolated curves in region 3, the amplitude of these periodic solutions as $z\to\pm\infty$ goes to zero, thus yielding truly localized solitary waves. These are known as embedded solitons \cite{cm99} and will be investigated for our system (2.1) in future work.
\item[c.] {\bf Near $C_2$}: In this region, where $\la_{1-4}=\pm i\omega$, $\pm i\omega$, derivation and analysis of a complicated normal form \cite{etbci87}-\cite{ip93} in Section 5 shows the possible occurrence of so-called \textit{envelope} homoclinic solutions of the form  sech$(kt)e^{(i\gamma\theta)}$ (and with oscillating tails) in the so-called \textit{subcritical} form. However, occurrence or persistence of these solutions in the full nonlinear system (2.1) is a non-trivial issue (and each system must be analyzed separately \cite{ip93}-\cite{ik90}). The persistence in (2.1) is considered in various ways in \cite{ip93}-\cite{ik90}. Indeed, there is no problem for solutions having \textit{one hump or peak}. The only open problem is for nonsymmetric solutions.
\item[d.] {\bf Near $C_3$:} There is no small-amplitude bifurcation on $C_3$, on which $\la_{1-4} =\pm\la , \pm\la$ and the fixed point \eqref{fixed_pt} remains hyperbolic. However, as we discuss below, there is a bifurcation across it causing the creation of an infinite multiplicity of homoclinic orbits.
\end{itemize}

We turn next to each of the regions 1 to 4 in Figure 1 to discuss the possible occurrence and multiplicity of homoclinic orbits in each one.
\begin{itemize}
\item[a.] {\bf Region 1}. In this region $\la_{1-4}=\pm\la\pm i\omega$ and fixed point \eqref{fixed_pt} is a saddle-focus. Using a Shil$'$nikov type analysis, one may show \cite{jh93},\cite{jh97} for general reversible systems such as (2.1) that the existence of one symmetric homoclinic orbit implies the existence of an infinity of others. Hence, we expect our system (2.1) to admit an infinity of such symmetric $N$-pulses for each $N>1$. Here, a symmetric $N$-pulse oscillates $N$ times in phase-space for $z\in (-\infty ,\infty )$ (or, more technically, crosses a transversal section to the primary symmetric 1 pulse $N$ times). In the context of the Ostrovsky equation, these would be $N$-peaked solitary waves, and we expect an infinite family for all $N>1$ for parameters $(p,q)$ in Region 1 of Figure 1.
\item[b.] {\bf Region 2}. In this region $\la_{1-4}=\pm\la_1,\pm\la_2$ and fixed point \eqref{fixed_pt} is a hyperbolic saddle point. Thus, there is no a priori reason for multiplicity of homoclinic orbits in this region. However, depending on the actual form of the nonlinear term, a symmetric homoclinic orbit to \eqref{fixed_pt} may exist (see \cite{etbci87},\cite{ip93}). Also, depending on further conditions \cite{bg96},\cite{ik90}, a further \textit{orbit-flip} bifurcation may cause complex dynamics in its neighborhood. In the context of our system (2.1), these issues will need further investigation to establish possible existence of solitary wave solutions in this region of $(b,a)$ space.
\item[c.] {\bf Region 3.} The generic situation in this region has already been considered in the discussion above pertaining to the region near curve $C_1$. As mentioned there, the structure and multiplicity of the delocalized solitons in Region 3, as well as the existence of embedded solitons on isolated curves, will be investigated for the Ostrovsky equation in future work.
\item[d.] {\bf Region 4.} In this region, $\la_{1-4}=\pm i\omega_1$, $\pm i\omega_2$ and \eqref{fixed_pt} is a focus. No homoclinic orbits are known to exist in general here, although complex dynamics may occur \cite{jh93},\cite{jh97}. Special results exist for $\omega_1 \approx \omega_2$.
\end{itemize}

Having outlined the possible families of orbits homoclinic to fixed point \eqref{fixed_pt} of \eqref{lin_eq},
we proceed next to derive convergent series solutions for solutions homoclinic to this fixed point in Region 1 discussed above.

\section{Reversible case: symmetric homoclinic solutions}\label{gen}

Here we shall derive solutions in the case where the nonlinear terms are chosen so that \eqref{eq_genC}
is reversible as discussed in the preceding section. However, the method of undetermined coefficients we employ for
the analytic solutions proves robust, and hence, in the following section, the series solutions
developed here are generalized to include nonlinear terms such that
\eqref{eq_genC} may be non-reversible, and hence to
cases more general than may be treated within the qualitative, reversible systems
bifurcation framework which has been discussed up to this point.

In this section, we choose the nonlinear function $f$ in \eqref{eq_genC} in such a way that the equation is reversible. So as to
include the traveling wave reductions of various, physically-relevant nonlinear PDEs discussed in \cite{C98}, we choose
the nonlinearity to be of the form:

\begin{equation}\label{rev}
f(u, u^{'}, u^{''}, u^{'''})=c u u^{''}+d u^{'2}+g u^2+h u^3.
\end{equation}

\noindent The resulting equation:

\begin{equation}\label{eq_rev}
u^{''''}-b u^{''}+a u=c u u^{''}+d u^{'2}+g u^2+h u^3,
\end{equation}

\noindent is reversible under the standard reversibility of classical mechanical systems:

\begin{equation}\label{rev_rule}
z \rightarrow -z \qquad (u, u^{'}, u^{''}, u^{'''}) \rightarrow (u, -u^{'}, u^{''}, -u^{'''}).
\end{equation}

\noindent Mathematically, this translates to solutions having even parity in $z$.

The equation \eqref{eq_rev} has the following three fixed points:

\begin{equation}\label{equi}
u=0,\quad u=\frac{-g\pm\sqrt{g^2+4ah}}{2h},\qquad h\neq 0
\end{equation}

Notice that if $h=0$ (i.e. if there is no $u^3$ term) there are only two fixed points $u=0$ and $u=a/g$.

Let us now proceed to construct in region $1$ of Figure 1 Shilnikov-type homoclinic orbits to the origin. Suppose that for $z>0$:

\begin{equation}\label{up_sym}
u^+(z)=\sum_{k=\,1}^{\infty} a_k e^{k\alpha z},
\end{equation}

\noindent where $\alpha<0$ is an undetermined constant and $a_k,\, k\geq 1$ are, at the outset, arbitrary coefficients.
Substituting the series \eqref{up_sym} into the reversible equation \eqref{eq_rev} yields:

\begin{equation}\label{sub_ser_sym}
\sum_{k=\,1}^{\infty}\left((k\alpha)^4-b(k\alpha)^2+a\right)a_ke^{k\alpha z}=F_1+F_2,
\end{equation}

\noindent with:

\begin{equation}
F_1=\sum_{k=\,2}^{\infty}F_1^{(k)},\qquad\qquad F_2=\sum_{k=\,3}^{\infty}F_2^{(k)}
\end{equation}

\noindent where

\begin{eqnarray}\label{F1F21}
F_1^{(k)}&=&\left\{\begin{array}{ll}0 \hskip 8.6cm  {\rm for\ }k=1\\
\,\\
\sum_{i=1}^{k-1}\left(c(k-i)^2\alpha^2+d(k-i)i\alpha^2+g\right)\,a_{k-i}\,a_i\, e^{k\alpha z}, \quad {\rm for\ }k>1
\end{array}
\right.\\\nonumber
\,\\\label{F1F22}
F_2^{(k)}&=&\left\{\begin{array}{ll}0 \hskip 5.55cm {\rm for\ }k=1,2\\
\,\\
\sum_{j=2}^{k-1}\sum_{l=1}^{j-1}h\, a_{k-j}\,a_{j-l}\,a_l\, e^{k\alpha z}\qquad {\rm for\ }k>2
\end{array}
\right.
\end{eqnarray}

\noindent Comparing the coefficients of $e^{k\alpha z}$ for each $k$, one has for $k=1$:

\begin{equation}\label{eq_roots}
\left(\alpha^4-b\, \alpha^2+a\right)a_1=0.
\end{equation}

\noindent Assuming $a_1\neq0$ (otherwise $a_k=0$ for all $k>1$ by induction), results in the values:

\begin{equation}\label{roots}
\alpha_{1-4}=\pm\sqrt{\frac{b\pm\sqrt{b^2-4a}}{2}},
\end{equation}

\noindent which, in region 1 of Fig.\ref{parabola}, are of the form $\alpha_1=-\lambda+i\omega, \alpha_2=\bar{\alpha}_1,
\alpha_3=\lambda+i\omega, \alpha_4=\bar{\alpha}_3$.

As our series solution \eqref{up_sym} needs to converge for $z>0$, we pick the two roots $\alpha_{1,2}$ with negative real parts:

\begin{equation}\label{up_symc}
u^+(z)=\sum_{k=\,1}^{\infty} a_k e^{k\alpha_1 z}+c.c.
\end{equation}

\noindent For $k>1$ one has:

\begin{equation}\label{rel_coef}
p(k\alpha_1)a_k=F_1^{(k)}+F_2^{(k)},
\end{equation}

\noindent where the polynomial:

\begin{equation}\label{pol}
p(k\alpha_1)=(k\alpha_1)^4-b(k\alpha_1)^2+a,
\end{equation}

\noindent is different from zero for $k>1$, as may be seen from \eqref{eq_roots}. Moreover, the terms into $F_1^{(k)}$ and $F_2^{(k)}$ involves series coefficients $a_j$, with $j>k$, as can be easily observed from the expressions \eqref{F1F21}-\eqref{F1F22}.

Therefore, the coefficients $a_k,\,k>1$ can be directly obtained from the equation \eqref{rel_coef} as follows:

\begin{equation}\label{coeff_rel}
a_k=\varphi_k\,a_1^k,
\end{equation}

\noindent where $\varphi_k,\, k>1$ are known functions given in Appendix \ref{app1} and
depending on $\alpha_1$ and the coefficients of the equation \eqref{eq_rev}.

The first part of the homoclinic orbit corresponding to $z>0$ has thus been determined in terms of $a_1$:

\begin{equation}\label{up_sym_comp}
u^+(z)=a_1 e^{\alpha_1 z}+\sum_{k=\,2}^{\infty}\varphi_k a_1^k e^{k\alpha_1 z}+c.c.
\end{equation}

\noindent We shall now construct the second part corresponding to $z<0$. We remark that since the equation \eqref{eq_rev} is reversible with respect to the transformation \eqref{rev_rule}, then for all $z<0$ we have that:

\begin{equation}\label{um_sym_comp}
u^-(z)=b_1 e^{\alpha_4 z}+\sum_{k=\,2}^{\infty}\varphi_k b_1^k e^{k\alpha_4 z}+c.c.
\end{equation}

As we seek continuous solutions to \eqref{eq_rev} at $z=0$, i.e. :

\begin{equation}\label{cont_sym}
u^+(0)=a_1 +\sum_{k=\,2}^{\infty}\varphi_k a_1^k=0=b_1 +\sum_{k=\,2}^{\infty}\varphi_k b_1^k=u^-(0)
\end{equation}

\noindent it is sufficient to choose $b_1=a_1$ and the relation

\begin{equation}\label{eq_a1}
a_1 +\sum_{k=\,2}^{\infty}\varphi_k a_1^k=0
\end{equation}

\noindent must hold. As the $\varphi_k$ are uniquely determined in terms of $\alpha_1$ and the coefficients of the equation \eqref{eq_rev}, the relation \eqref{eq_a1} allows one to compute $a_1$ (in practice this must be done numerically at fixed values of the coefficients of the equation \eqref{eq_rev}).

We have thus determined that the equation \eqref{eq_rev} has a homoclinic orbit to the origin of the form:

\begin{equation}\label{orbit}
u(z)=\left\{\begin{array}{lll}
u^+(z)\qquad z>0\\
0\qquad\qquad z=0\\
u^-(z)\qquad z<0
\end{array}\right.
\end{equation}

However, the converegence of the above series is still unclear. So as to not get mired in details
and lose the main thread of our development, we defer proof of convergence to Section 6,
assuming its validity for the present. Note also that the above solution in \eqref{orbit}
represents homoclinic orbits at ALL $(b,a)$ points in Region 1 of Figure 1.

\begin{figure}[h]
\begin{center}
\subfigure[] {\epsfxsize=2.6 in \epsfbox{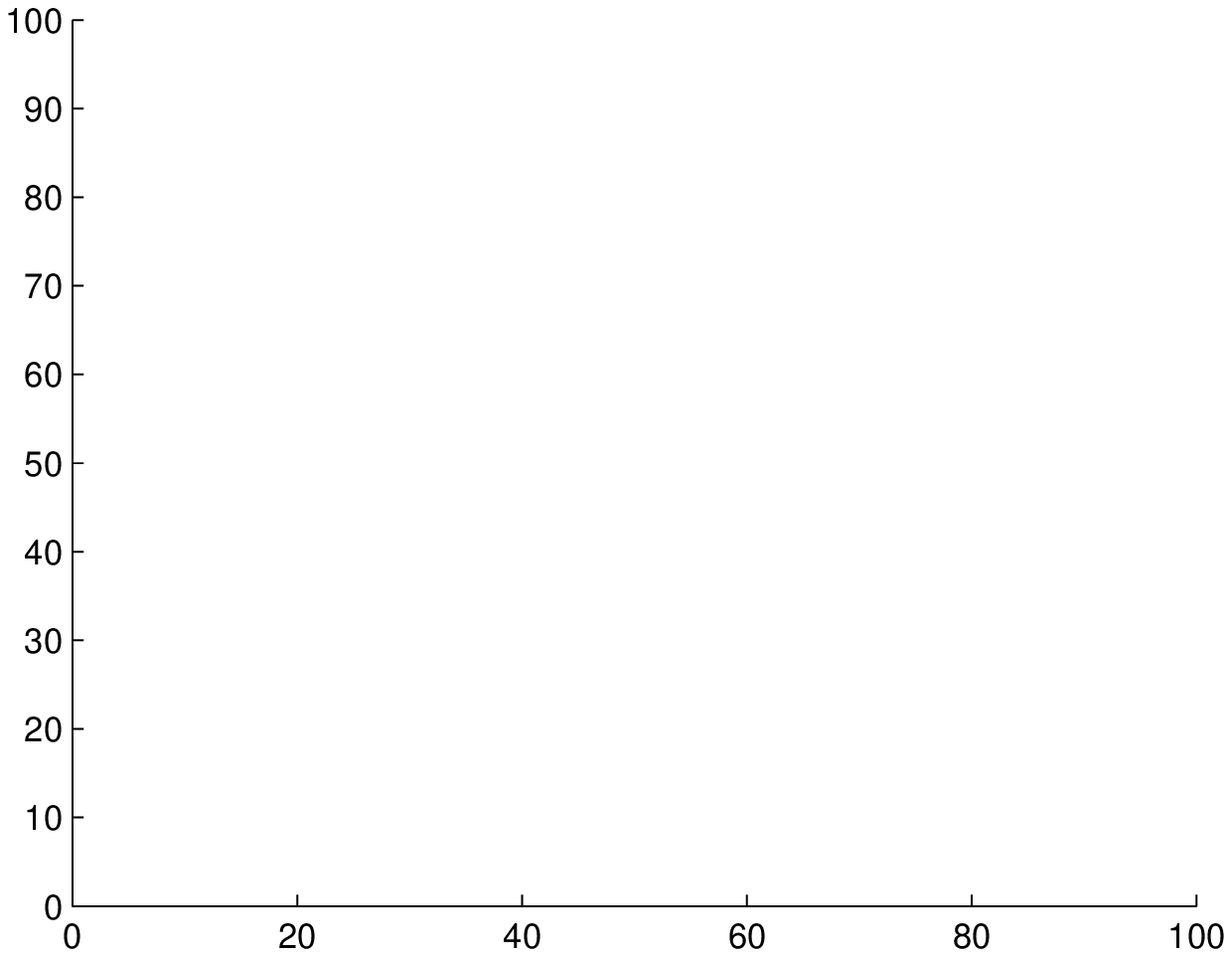}}
\subfigure[] {\epsfxsize=2.6
 in \epsfbox{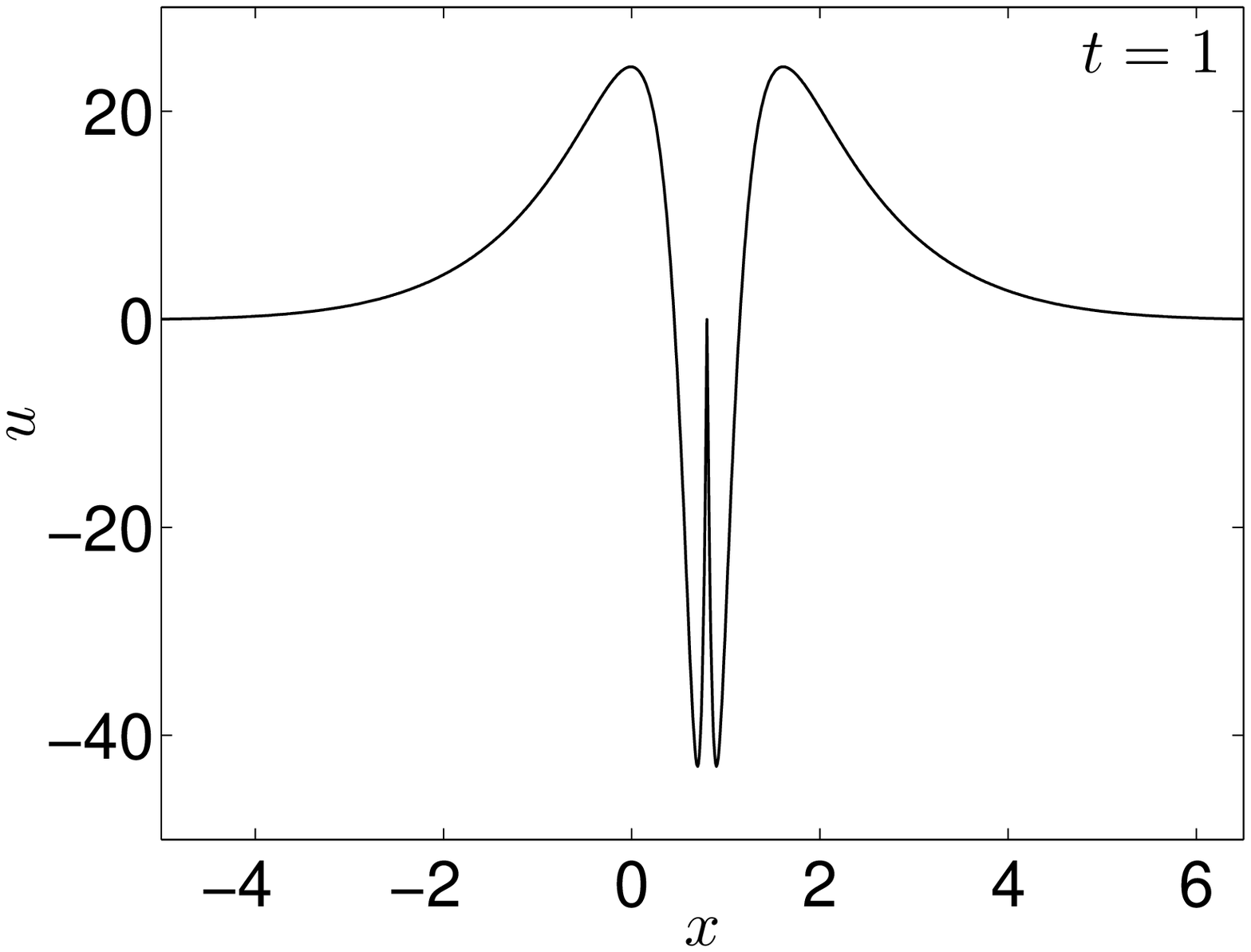}}
\end{center}
\caption{\label{1sym} Symmetric homoclinic solutions. For the parameter choice see the text.}
\end{figure}

In Fig.\ref{1sym}-(a) we fix $a=0.8$ and $b=1.5$ in region 1 and the other parameters are $c=0.2, d=0.1, g=0.05$ and $h=0.02$. The series coefficient $a_1$ has been chosen as the root $a_1=40.4440 -14.2061i$ of \eqref{eq_a1}.
In Fig.\ref{1sym}-(b) all the equation parameters are chosen as in Fig.\ref{1sym}. The series coefficient $a_1$ has been chosen as the root $a_1=41.6756 - 0.9564i$ of \eqref{eq_a1}. In both figures the solution $u$ is shown with respect to the variable $x$, for fixed $t$.

\begin{figure}[h]
\begin{center}
\subfigure[] {\epsfxsize=2.6 in \epsfbox{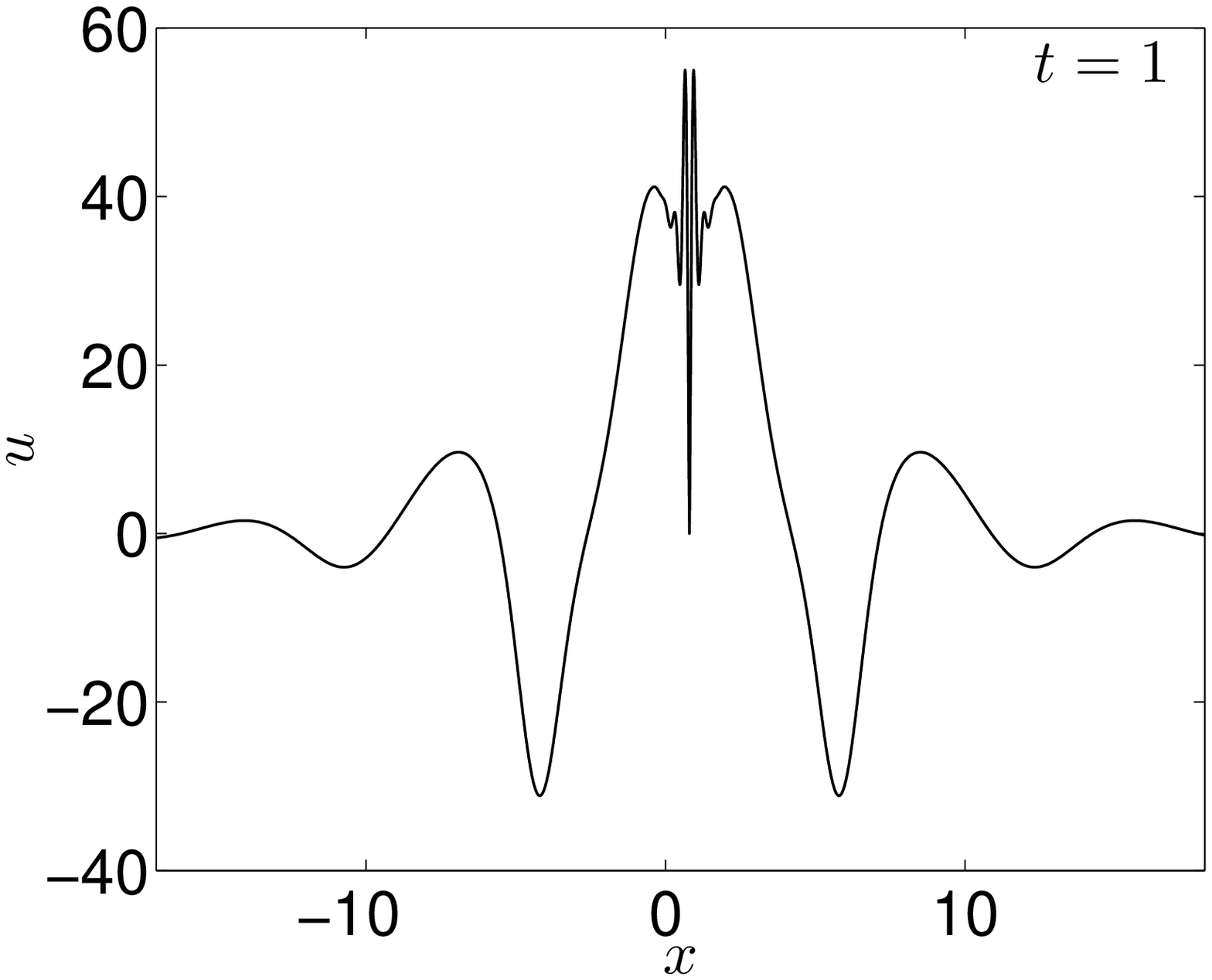}}
\subfigure[] {\epsfxsize=2.6
 in \epsfbox{vuoto.eps}}
\end{center}
\caption{\label{3sym} Symmetric homoclinic solutions: the pulses show a jagged part. For the parameter choice see the text.}
\end{figure}

In Fig.\ref{3sym}-(a) we fix $a=0.8$ and $b=-1.5$ in region 1 and the other parameters are as in Fig.\ref{1sym}. The series coefficient $a_1$ has been chosen as the root $a_1=14.814-41.6433i$ of \eqref{eq_a1}.
In Fig.\ref{3sym}-(b) all the equation parameters are chosen as in Fig.\ref{1sym}. The series coefficient $a_1$ has been chosen as the root $a_1=-12.0892 -44.4016i$ of \eqref{eq_a1}. In both figures the solution $u$ is again shown with respect to the variable $x$ at fixed $t$.

\begin{figure}[h]
\begin{center}
\subfigure[] {\epsfxsize=2.6 in \epsfbox{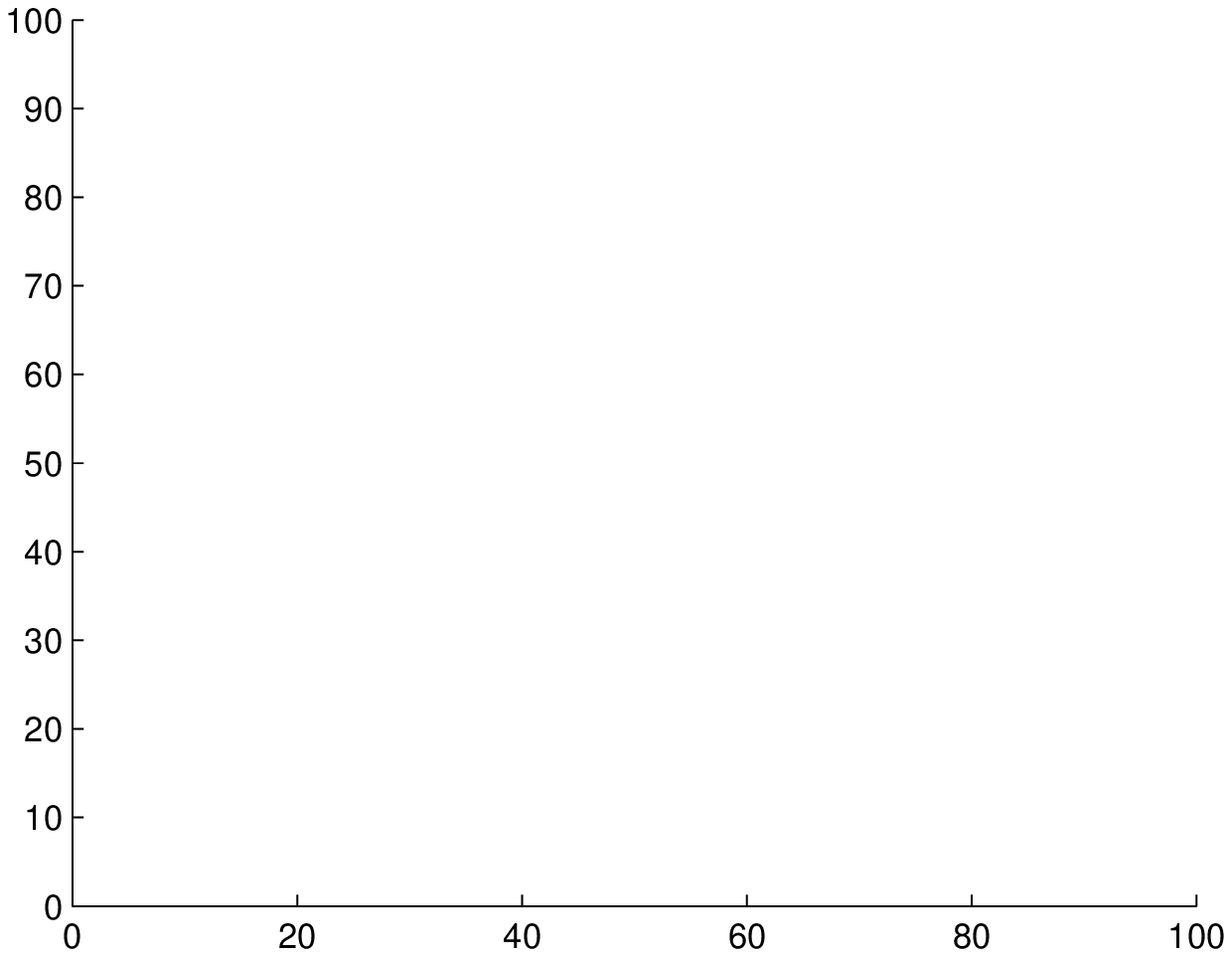}}
\subfigure[] {\epsfxsize=2.6
 in \epsfbox{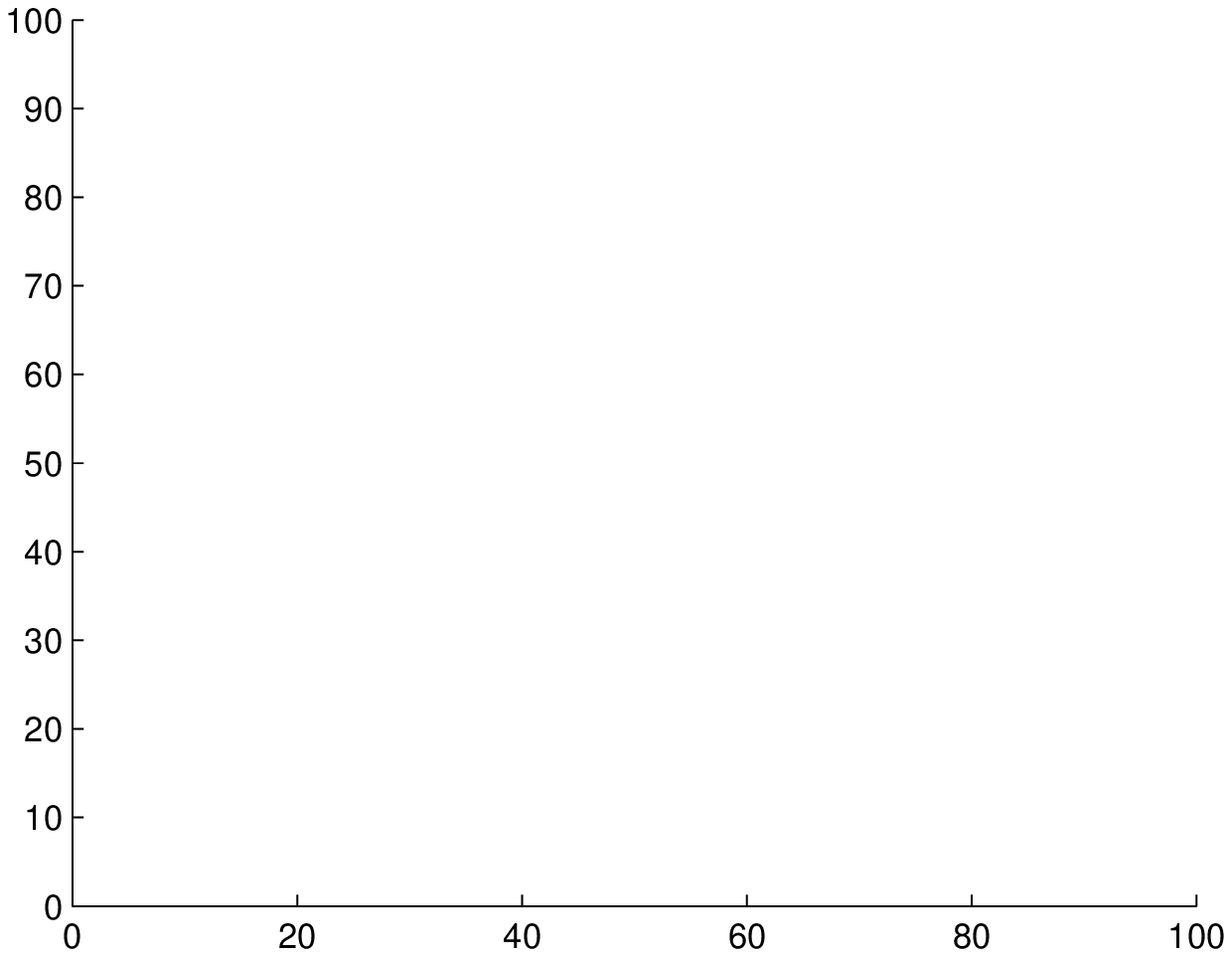}}
\end{center}
\caption{\label{4sym} Symmetric homoclinic solutions: the jagged part of the pulses becomes smooth choosing stronger nonlinearities.}
\end{figure}

Notice that the jagged part of the pulses in Fig.\ref{3sym} becomes smooth by changing the strength of the nonlinearities, as shown in Fig. \ref{4sym}, where $a$ and $b$ are still the same as in Fig.\ref{3sym}, but the other parameters are chosen as $c=0.7$, $d=0.6$, $g=0.8$ and $h=0.4$.

\begin{figure}[h]
\begin{center}
\subfigure[] {\epsfxsize=2.5 in \epsfbox{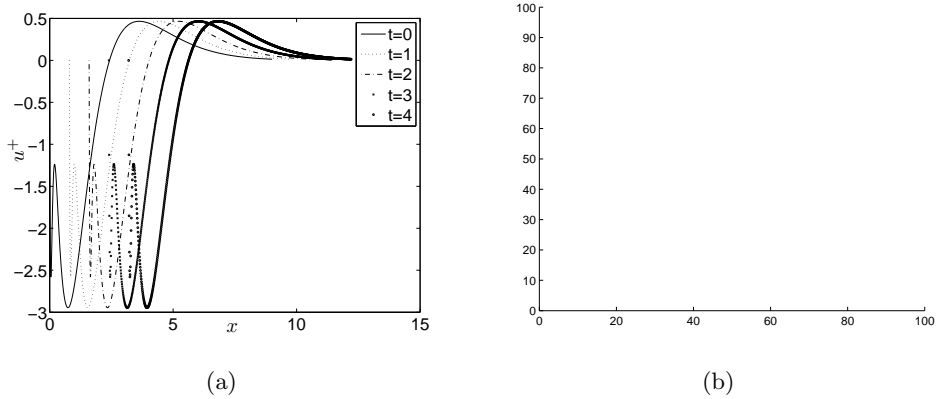}}
\subfigure[] {\epsfxsize=2.6
 in \epsfbox{vuoto.eps}}
\end{center}
\caption{\label{travel_sym} Traveling wave solution}
\end{figure}

\noindent In order to demonstrate the traveling-wave nature of our series solutions, Figure \ref{travel_sym} plots
the solutions at constant time $t=0,\dots,4$. In Fig.\ref{travel_sym}-(a) we fix $a=0.8$ and $b=1.5$ in region 1 and the other parameters are $c=0.3, d=0.1, g=0.05, h=0.02$. The series coefficient $a_1$ has been chosen as the root $a_1=-18.4550+15.8744i$ of \eqref{eq_a1}. In Fig.\ref{travel_sym}-(b) $a=3$ and $b=-1.5$ and the series coefficient $a_1=12.5579+36.8257i$, the other parameters are as in Fig.\ref{travel_sym}-(a).

In addition to the convergence of our multi-infinite series solutions which, as mentioned earlier, is detailed in Section 6, we may also check our results versus earlier results obtained by other methods,
for instance some derived numerically by homoclinic shooting \cite{Champ-Groves}. Scaling the governing equation in
\cite{Champ-Groves}, it is easy to see that only Figures 3 and 22 in that paper lie in Region 1 of our Figure 1.
The homoclinic orbit given by our series solutions \eqref{orbit} for parameters corresponding to those two figures in
\cite{Champ-Groves} are shown in Fig.\ref{e1} below, and they agree well with the numerical solutions in the earlier paper.

\begin{figure}[h]
\begin{center}
\subfigure[] {\epsfxsize=2.6 in \epsfbox{vuoto.eps}}
\subfigure[] {\epsfxsize=2.6
 in \epsfbox{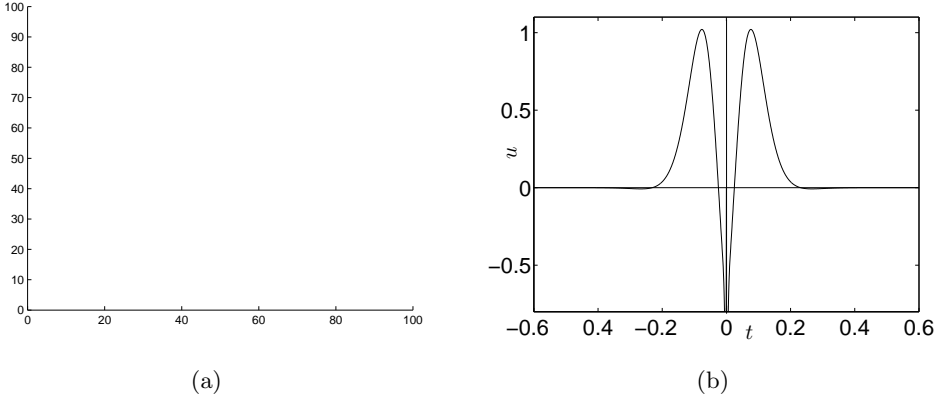}}
\end{center}
\caption{\label{e1} Comparison of our series solutions with the numerical solutions in \cite{Champ-Groves}. In (a) the parameters are $a=7.5, b, d=-3.75, c=-7.5, g=-11.25, h=0$. In (b) the parameters are $a=15, b=3.75, c=-7.5, d=-3.75,  g=-11.25, h=0$.}
\end{figure}

\section{Non-reversible case: asymmetric homoclinic solutions}

Since the series solutions in the previous section agree closely with numerical results for
the reversible case,
we now include other nonlinear terms in \eqref{eq_genC} which break the reversibility of the equation:

\begin{equation}\label{f_gen}
f(u, u^{'}, u^{''}, u^{'''})=c u u^{''}+d u^{'2}+g u^2+h u^3+p u u^{'''}+q u u^{''}+r u u^{'}+su^3u^{'}.
\end{equation}

\noindent The resulting equation is:

\begin{equation}\label{eq_gen}
u^{''''}-b u^{''}+a u=c u u^{''}+d u^{'2}+g u^2+h u^3+p u u^{'''}+q u u^{''}+r u u^{'}+su^3u^{'}.
\end{equation}

\noindent Let us build Shilnikov-type homoclinic orbits to the origin in region $1$ of Figure 1, as in the previous section. Suppose that, for $z>0$, the solution has the form \eqref{up_sym}.

Let us substitute the series \eqref{up_sym} into the equation \eqref{eq_gen} and it reduces to:

\begin{equation}\label{sub_ser}
\sum_{k=\,1}^{\infty}\left((k\alpha)^4-b(k\alpha)^2+a\right)a_ke^{k\alpha z}=F_1+F_2+F_3,
\end{equation}

\noindent with:

\begin{equation}
F_1=\sum_{k=\,2}^{\infty}F_1^{(k)},\qquad\qquad F_2=\sum_{k=\,3}^{\infty}F_2^{(k)},\qquad\qquad F_3=\sum_{k=\,4}^{\infty}F_3^{(k)}
\end{equation}

\noindent where:

\begin{eqnarray}\label{F1F2F31}
F_1^{(k)}&=&\left\{\begin{array}{ll}0 \hskip 8.6cm  {\rm for\ }k=1\\
\,\\
\sum_{i=1}^{k-1}\left(p(k-i)^3\alpha^3+(k-i)^2(c+iq\alpha)\alpha^2+\right.\hskip 1.6cm  {\rm for\ }k>1\\
\left.\hskip1.1cm r(k-i)\alpha+d(k-i)i\alpha^2+g\right)\,a_{k-i}\,a_i\, e^{k\alpha z}
\end{array}\right.\\\label{F1F2F32}
F_2^{(k)}&=&\left\{\begin{array}{ll}0 \hskip 8.6cm  {\rm for\ }k=1,2\\
\,\\
\sum_{j=2}^{k-1}\sum_{l=1}^{j-1}h\, a_{k-j}\,a_{j-l}\,a_l\, e^{k\alpha z}\hskip 3.8cm  {\rm for\ }k>1
\end{array}\right.\\\label{F1F2F33}
F_3^{(k)}&=&\left\{\begin{array}{ll}0 \hskip 8.6cm  {\rm for\ }k=1,2,3\\
\,\\
\sum_{i=3}^{k-1}\sum_{j=2}^{i-1}\sum_{l=1}^{j-1}s\, a_{k-i}\,a_{j-l}\,a_l\,a_{i-j}\, e^{k\alpha z}\hskip 2.2cm  {\rm for\ }k>1
\end{array}\right.
\end{eqnarray}

\noindent Comparing the coefficients of $e^{k\alpha z}$ for each $k$, one has for $k=1$ the equation \eqref{eq_roots}.
Once again we assume $a_1\neq0$ (otherwise $a_k=0$ by induction) for all $k>1$, resulting in the four complex roots \eqref{roots}.

As our series solution \eqref{up_sym} needs to converge for $z>0$, we pick the two roots $\alpha_{1,2}$ with negative real parts and the form of the solution for $z>0$ is still \eqref{up_symc}.
Moreover, for $k>1$ one obtains the form of the coefficients $a_k$ as in \eqref{coeff_rel},
where $\varphi_k,\, k>1$ are known functions depending on $\alpha_1$ and the coefficients of the equation \eqref{eq_gen}. The details are different now and are given in Appendix \ref{app2}.

The first part of the homoclinic orbit corresponding to $z>0$ has thus been determined in terms of $a_1$, formally as in \eqref{up_sym_comp}.

We shall now construct the second part corresponding to $z<0$. Since the equation is not reversible we do not have any symmetry property for the solution, therefore impose the following solution form:

\begin{equation}\label{um}
u^-(z)=\sum_{k=\,1}^{\infty}b_k e^{k\alpha z}+c.c.
\end{equation}

\noindent where $\mathfrak{Re}(\alpha)>0$, because the solution $u^-(z)$ needs to converge for $z<0$.
Working as for $z>0$, we obtain for $k=1$ the following equation:

\begin{equation}\label{eq_rootsB}
\left(\alpha^4-b\, \alpha^2+a\right)b_1=0.
\end{equation}

\noindent Assuming $b_1\neq0$ (otherwise $b_k=0$ for all $k>1$ by induction), we find the four roots in \eqref{roots} and, for $z>0$, we choose $\alpha=\alpha_4$.
Moreover, for $k>1$ one obtains the following equation:

\begin{equation}\label{rel_coefb}
p(k\alpha_4)b_k=F_1^{(k)}+F_2^{(k)}+F_3^{(k)},
\end{equation}

\noindent where the polynomial

\begin{equation}\label{pol_b}
p(k\alpha_4)=(k\alpha)^4-b(k\alpha)^2+a,
\end{equation}

\noindent is non-zero for $k>1$ and the quantities $F_1^{(k)}, F_2^{(k)}$ and $F_3^{(k)}$ are formally the same as in \eqref{F1F2F31}-\eqref{F1F2F33} once the coefficients $a_k$ are substituted by $b_k$.

Therefore, the series coefficients can be easily obtained as follows:

\begin{equation}\label{coeff_relB}
b_k=\psi_k\,b_1^k,
\end{equation}

\noindent where $\psi_k$ are given in terms of $b_1$ and the coefficients of equations \eqref{eq_gen}. Notice that in this case there is no relation between $\varphi_k$ and $\psi_k$ as in the reversible case (see the details in Appendix \ref{app2}).

 As we want to construct a solution of the following form:

\begin{equation}\label{solutionB}
u(z)=\left\{\begin{array}{lll}
u^+(z)\qquad z>0\\
0\qquad\qquad z=0\\
u^-(z)\qquad z<0
\end{array}\right.
\end{equation}

\noindent which is continuous at $z=0$, we impose:

\begin{eqnarray}\label{contB1}
u^+(0)&=&a_1 +\sum_{k=\,2}^{\infty}\varphi_k a_1^k=0\\\label{contB2}
u^-(0)&=&b_1 +\sum_{k=\,2}^{\infty}\psi_k b_1^k=0.
\end{eqnarray}

Hence, we choose $a_1$ and $b_1$ as the nontrivial solutions of the following polynomial equations:

\begin{eqnarray}\label{eq_a1gen}
a_1 +\sum_{k=\,2}^{\infty}\varphi_k a_1^k&=&0,\\\label{eq_b1gen}
b_1 +\sum_{k=\,2}^{\infty}\psi_k b_1^k&=&0.
\end{eqnarray}

\noindent In practice the equations \eqref{eq_a1gen}-\eqref{eq_b1gen} are numerically solved at fixed values of the coefficients of the equation \eqref{eq_gen}. These solutions are not unique. In the following we show some numerical simulations of the homoclinic orbit \eqref{solutionB} for different values of the parameters of the equation \eqref{eq_gen}.

\begin{figure}[h]
\begin{center}
\subfigure[] {\epsfxsize=2.6 in \epsfbox{vuoto.eps}}
\subfigure[] {\epsfxsize=2.6
 in \epsfbox{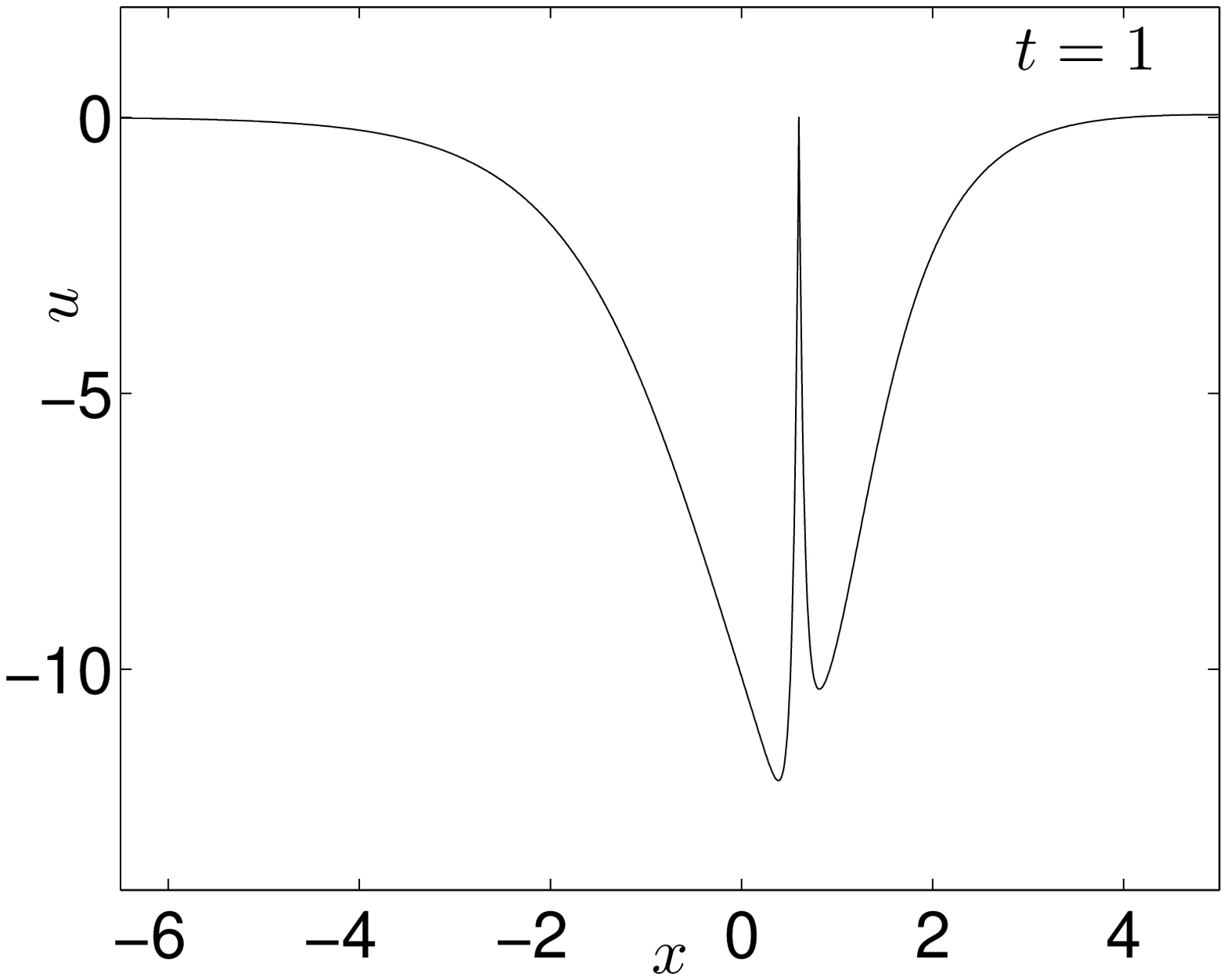}}
\end{center}
\caption{\label{1gen} Asymmetric homoclinic solutions. For the parameter choice see the text.}
\end{figure}

In Fig.\ref{1gen} the parameters $a=0.6$ and $b=1.5$ are chosen in Region 1 of Fig.\ref{parabola}. The other parameters are fixed as $p=0.1$, $q=0.3$, $c=0.2$, $r=0.1$, $d=0.3$, $g=0.5$, $h=0.2$ and $s=0.1$. In Fig.\ref{1gen}-(a) the solutions of equations \eqref{eq_a1gen} and \eqref{eq_b1gen} are respectively chosen as $a_1=2.3840 - 8.9933i$ and $b_1=11.8609 + 3.8370i$; in Fig.\ref{1gen}-(b) they are chosen as $a_1= -3.1533 - 8.6762i$ and $b_1=-10.4363 - 6.2990i$.

\begin{figure}[h]
\begin{center}
\subfigure[] {\epsfxsize=2.6 in \epsfbox{vuoto.eps}}
\subfigure[] {\epsfxsize=2.6
 in \epsfbox{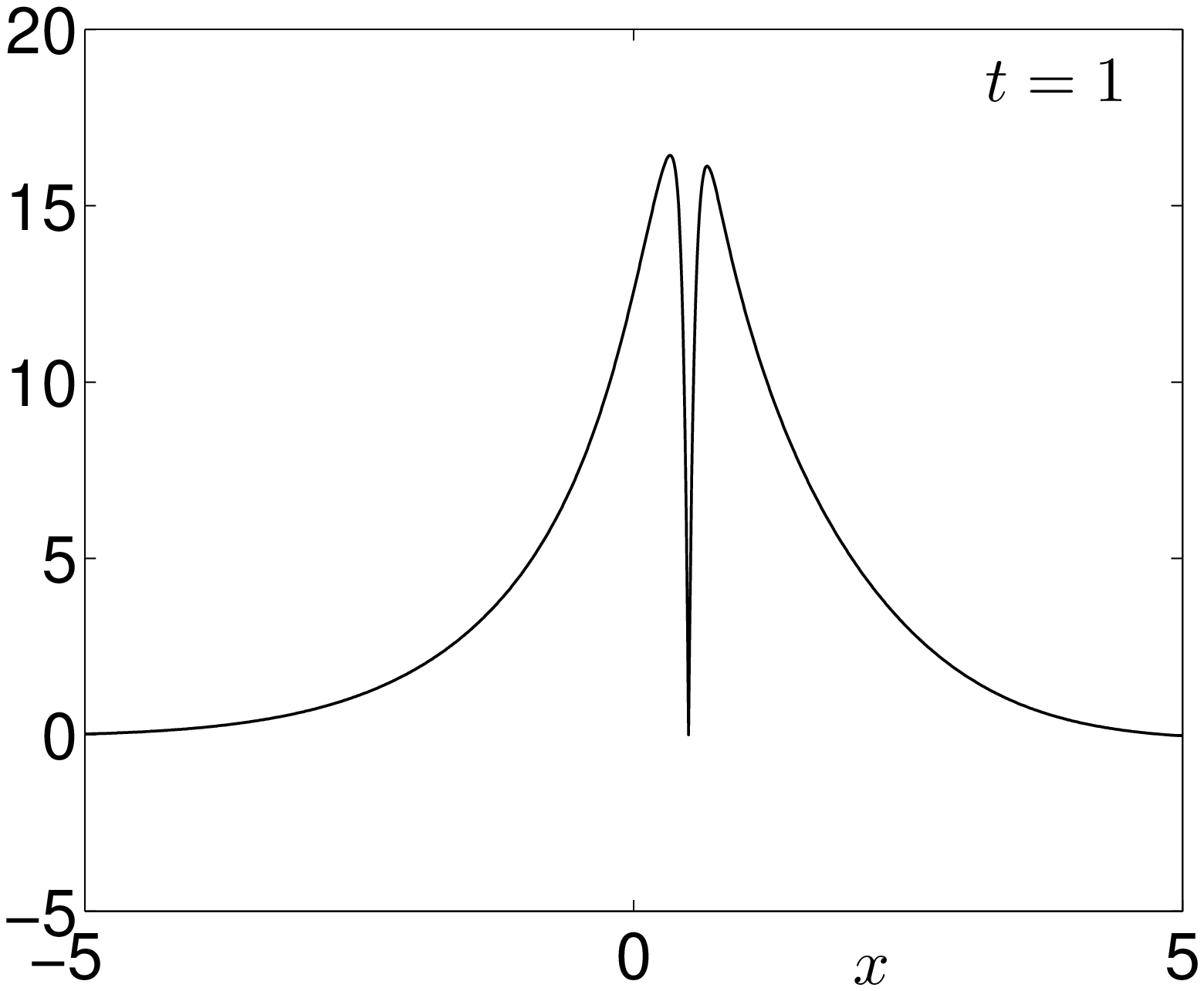}}
\end{center}
\caption{\label{2gen} Asymmetric homoclinic solutions. For the parameter choice see the text.}
\end{figure}

In Fig.\ref{2gen} the parameters $a=0.5$ and $b=1$ are chosen into the region 1 of Fig.\ref{parabola}. The other parameters are fixed as $p=0.01$, $q=0.3$, $c=0.02$, $r=0.1$, $d=0.3$, $g=0.05$, $h=0.2$ and $s=0.1$. In Fig.\ref{1gen}-(a) the solutions of equations \eqref{eq_a1gen} and \eqref{eq_b1gen} are respectively chosen as $a_1=-7.6505 - 6.8582i$ and $b_1=4.7888 + 6.1539i$; in Fig.\ref{2gen}-(b) they are chosen as $a_1=10.7954 + 1.8878i$ and $b_1=7.4002 - 2.3540i$.

\begin{figure}[h]
\begin{center}
\subfigure[] {\epsfxsize=2.6 in \epsfbox{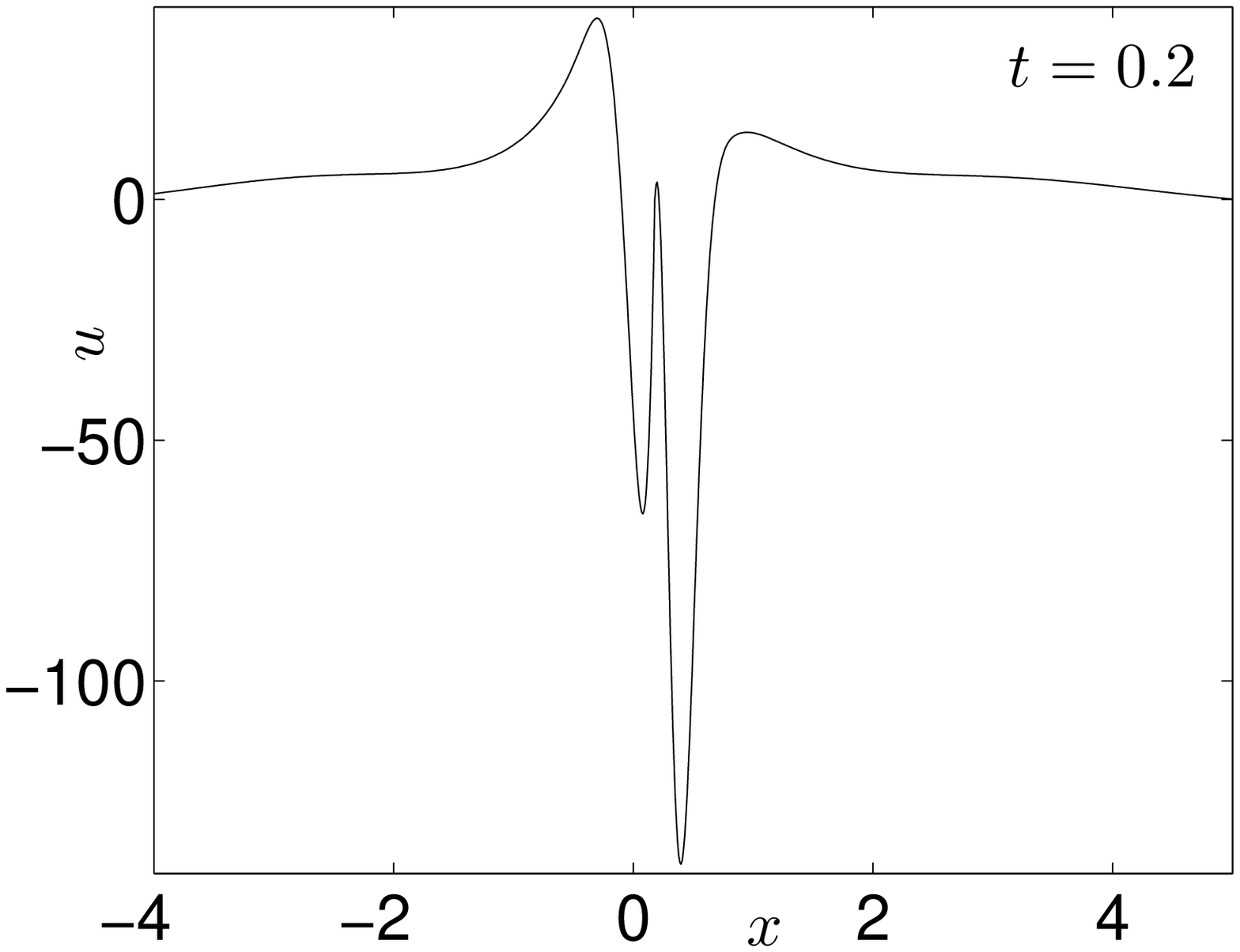}}
\subfigure[] {\epsfxsize=2.6
 in \epsfbox{vuoto.eps}}
\end{center}
\caption{\label{3gen} Asymmetric homoclinic solutions. For the parameter choice see the text.}
\end{figure}

In Fig.\ref{3gen} the parameters $a=0.9$ and $b=-1.2$ are chosen in Region 1 of Fig.\ref{parabola}. The other parameters are fixed as
$p=0.01$, $q=0.02$, $c=0.7$, $r=0.1$, $d=0.6$, $g=0.8$, $h=0.5$ and $s=0.03$. In Fig.\ref{1gen}-(a) the solutions of equations \eqref{eq_a1gen} and \eqref{eq_b1gen} are respectively chosen as $a_1=-8.5158 - 4.7092i$ and $b_1=-7.4718 - 7.3032i$; in Fig.\ref{3gen}-(b) they are chosen as $a_1=-8.5158 - 4.7092i$ and $b_1=-6.7027 -11.8976i$.

\begin{figure}[h]
\begin{center}
\subfigure[] {\epsfxsize=2.6 in \epsfbox{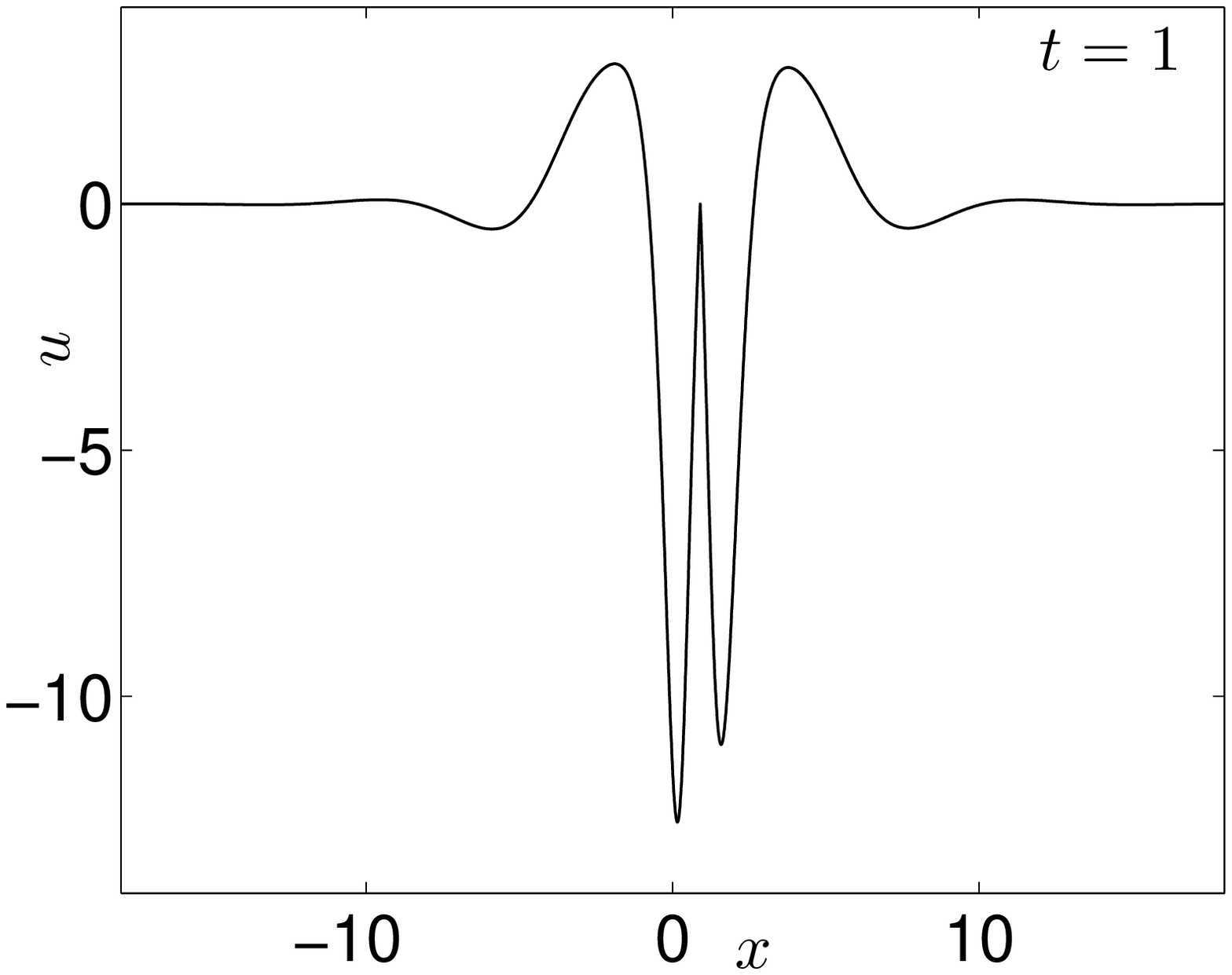}}
\subfigure[] {\epsfxsize=2.6
 in \epsfbox{vuoto.eps}}
\end{center}
\caption{\label{4gen} Asymmetric homoclinic solutions. For the parameter choice see the text.}
\end{figure}

Further examples are shown in Fig.\ref{4gen}. In particular, in Fig\ref{4gen}-(a) the parameters are chosen as $a=0.9$ and
$b=-1.0$ in region 1 of Fig.\ref{parabola}; and $p=0.7$, $q=0.6$, $c=0.8$,
$r=0.5$, $d=0.01$, $g=-0.02$, $h=0.1$ and $s=-0.03$. The roots of equations \eqref{eq_a1gen}-\eqref{eq_b1gen} are $a_1=-6.9709$ and $b_1=-7.2113$.
In Fig\ref{4gen}-(b) the parameters are chosen as $a=0.9$ and
$b=-1.6$ in region 1 represented in Fig.\ref{parabola}; and $p=0.8$, $q=0.6$, $c=0.85$,
$r=0.7$, $d=0.01$, $g=-0.02$, $h=0.1$ and $s=-0.03$. The roots of equations \eqref{eq_a1gen}-\eqref{eq_b1gen} are $a_1=-6.9709$ and $b_1=-7.2113$.

\begin{figure}[h]
\begin{center}
\subfigure[] {\epsfxsize=2.3 in \epsfbox{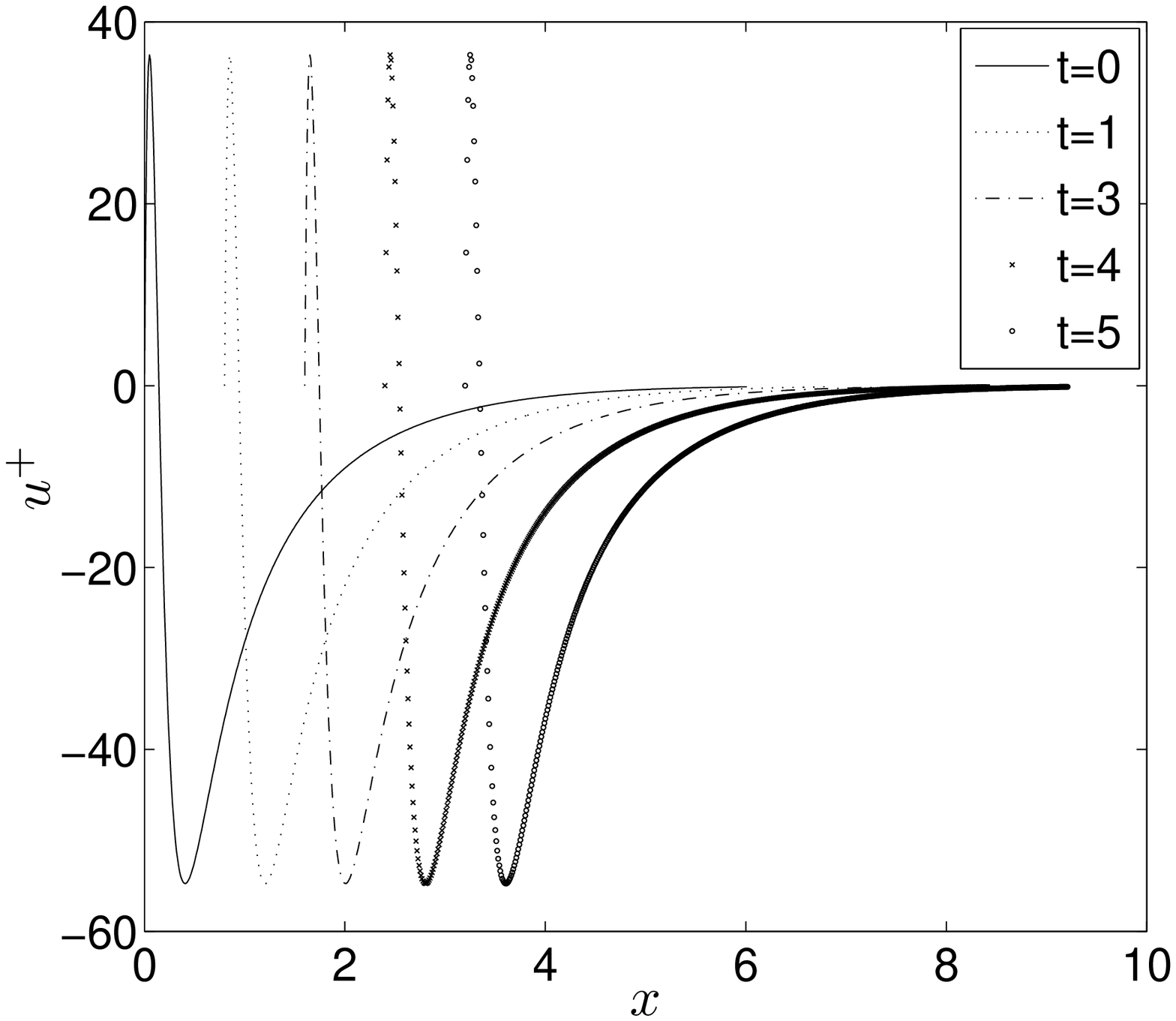}}
\subfigure[] {\epsfysize=2.0
 in \epsfbox{vuoto.eps}}
\end{center}
\caption{\label{trav1gen} Traveling wave solution}
\end{figure}

In order to demonstrate the traveling-wave nature of our series solutions as in the previous section, Figure \ref{trav1gen}  plots
the solutions at constant time $t=0,\dots,4$. Notice that we are only plotting the right solutions $u^+$ to avoid mess into the figure. The parameters are chosen as $c=0.3, d=0.1, g=0.05, h=0.02, p=0.1, q=0.2, r=0.1, s=0$; in Fig.\ref{trav1gen}-(a) we fix $a=0.8$ and $b=1.5$  in region 1 and the series coefficient $a_1=-16.0853-18.0011i$, in Fig.\ref{trav1gen}-(b) $a=0.8$ and $b=-1.5$ and $a_1=4.8877+29.9708i$.

\section{Proof of convergence}

We now return to the part of our derivation which had been deferred, viz. the proof that the series solutions
obtained for the homoclinic orbits in the two previous sections are indeed convergent. In particular, as part of this
proof, we will demonstrate the feature which was observed numerically earlier, i.e. that the later coefficients in our
multi-infinite-series solutions drop off rapidly.

As discussed in the previous two sections (see \eqref{coeff_rel} and \eqref{coeff_relB} respectively), one may inductively obtain
the $k$-th coefficient in the series in the form

\begin{equation}
a_k=\varphi_k\,a_1^k.
\end{equation}

More specifically (here we report the expressions into the symmetric case), for $k=2$:

\begin{equation}\label{step1}
a_2=\frac{\left((c+d)\alpha_1^2+g\right)\,a_1^2}{p(2\alpha_1)},
\end{equation}

\noindent For $k=3$:

\begin{equation}\label{step2}
a_3=\frac{\left((5c+4d)\alpha_1^2+2g\right)\,a_1a_2+h\,a_1^3}{p(3\alpha_1)}.
\end{equation}

\noindent For $k=4$:

\begin{equation}\label{step3}
a_4= \frac{(10c+ 6d)\alpha_1^2+2g)\,a_1a_3+(4(c+d) \alpha_1^2 +  g)\,a_2^2+ 3ha_1^2a_2}{p(4\alpha_1)}
\end{equation}

\noindent For $k=5$:

\begin{equation}\label{step4}
a_5= \frac{((17c + 8d)\alpha_1^2 + 2g)\,a_1a_4+((13c+ 12d)\alpha_1^2 + 2g)\,a_2a_3
+ 3ha_1a_3 + 3ha_1a_2}{p(5\alpha_1)}
\end{equation}

\noindent And so on for $k>5$:

\begin{equation}\label{ak}
a_k=\frac{F_1^{(k)}+F_2^{(k)}}{p(k\alpha_1)},
\end{equation}

\noindent where $F_1^{(k)}$ and $F_2^{(k)}$ can be inductively obtained as the product of $a_1^k$ with known functions of $\alpha_1$ and the parameters of the equation.

Following some standard steps in \cite{VC}-\cite{ZCT}, it is straightforward to obtain the bound

\begin{equation}
|a_k| < r(c,d,\xi)^{-(k+1)}  |a_1|^k,
\end{equation}

\noindent
where $r(c,d,\xi)$ is a constant dependent on the parameters $c$ and $d$ of our ODE, as well as Euler's constant $\xi$.
Note that this validates the steep fall-off in the magnitude of the late terms of our series solutions
which was numerically observed in Sections 4 and 5.

For typical system parameters considered in Sections 4 and 5, the above bound may
be expressed (evaluated) as:

\begin{equation}
|a_k| < 10^{-(k+1)}  |a_1|^k, k > 4.
\end{equation}

Hence, $|a_k|, k>4$ are bounded by some constant $L>0$. Using this, we may bound the series of absolute values associated to
the two halves  \eqref{up_sym} and \eqref{um_sym_comp} of our homoclinic orbit
solutions by $L$ times a convergent exponential (geometric) series. Hence our
series solutions are absolutely convergent, and hence convergent, by the \textit{Comparison Test}.

As a final comment, the fixed point $(b,a) = (0,0)$ of our ODE system
is, as detailed earlier in Section 3, a saddle-focus in Region 1 of Figure 1. Hence we may
conclude that our series solutions  \eqref{up_sym} and \eqref{um_sym_comp} represent the two halves
of the homoclinic orbits for all $(b,a)$ values in that region.

\section{Conclusions}

In summary, we have derived convergent, multi-infinite, series solutions for the homoclinic orbits
of our canonical fourth-order ODE system in Region 1 of Figure 1, for both reversible and non-reversible cases. If the ODE is a traveling-wave reduction of some PDE, these solutions would correspond to genuine or localized pulses of the PDE.
Given the many systems to which our canonical ODE pertains, our results are clearly of topical interest. In addition, the results for the non-reversible case go beyond those that have been typically
considered in most analyses conducted within bifurcation-theoretic settings.

As discussed in Section 3, the homoclinic solutions in Region 3 of Figure 1 are
known to be homoclinic or asymptotic to
exponentially small periodic orbits (corresponding to \textit{delocalized} solitary waves with exponentially small tails of
any corresponding PDE(s)). Considering such solutions is clearly a natural next step for future work
along the lines followed here for Region 1 of Figure 1. One possible approach
might be within the framework of \cite{JY01}, treating our series solutions here as
the unperturbed solutions $u_0$ in that work, and perturbing around them to derive
the tail amplitudes of delocalized solitary waves in Region 3 of our Figure 1.
Setting the tail amplitudes to zero would then yield the isolated curves
in $(b,a)$ parameter space (so far only known numerically) on which these delocalized waves become genuine \textit{embedded solitons}
\cite{cm99}. However,
the steps in the calculations would now involve a \textit{virtuoso} use of multi-infinite series, and we defer this to a future paper.

In another forthcoming paper the solutions and results here will be juxtaposed with, and placed in the framework of, a comprehensive categorization of ALL possible traveling-wave solutions, both smooth and non-smooth, for our governing canonical ODE.

\appendix
\section{Details of computations}\label{AppA}
\setcounter{equation}{0}
\subsection{Details for the symmetric case}\label{app1}

%
%
%
%
%
%
%
%
%
%
%

\subsection{Details for the general case}\label{app2}

For $k=2$ the equation \eqref{rel_coefb} reduces to:

\begin{equation}\label{step1g}
p(2\alpha_1)a_2=\left((p+q)\alpha_1^3+(c+d)\alpha_1^2+r\alpha_1+g\right)\,a_1^2,
\end{equation}

\noindent which gives $a_2$ in terms of $a_1$, in the form given in \eqref{coeff_rel}, where:

\begin{equation}\label{phi2g}
\varphi_2=\frac{(p+q)\alpha_1^3+(c+d)\alpha_1^2+r\alpha_1+g}{p(2\alpha_1)}.
\end{equation}

Analogously, for $k=3$ the equation \eqref{rel_coefb} reduces to:

\begin{equation}\label{step2g}
p(3\alpha_1)a_3=\left((9p+6q)\alpha_1^3+(5c+4d)\alpha_1^2+3r\alpha_1+2g\right)\,a_1a_2+h\,a_1^3.
\end{equation}

\noindent Taking into account the expression for $a_2$ in terms of $a_1$ obtained via \eqref{step1g}, the equation \eqref{step2g} gives $a_3$ in terms of $a_1$, in the form given in \eqref{coeff_rel}, where:

\begin{eqnarray}\label{phi3g}\nonumber
\varphi_3&=&\frac{1}{p(3\alpha_1)}\left[3(3p^2+5pq+2q^2)\alpha_1^6+(c(14p+11q)+d(13p+10q))\alpha_1^5+\right.\\
&\ &(3r(4p+3q)+c(5c+9d)+4d^2+16h)\alpha_1^4+(g(11p+8q)+\\\nonumber
&\ &\left.r(8c+7d))\alpha_1^3+
(g(7c+6d)-4bh+3r^2)\alpha_1^2+5gr\alpha_1+2g^2+ah\right]
\end{eqnarray}

\noindent And so on for $k>3$:

\begin{equation}\label{akg}
a_k=\frac{F_1^{(k)}+F_2^{(k)}+F_3^{(k)}}{p(k\alpha_1)},
\end{equation}

\noindent where $F_1^{(k)}$, $F_2^{(k)}$ and $F_3^{(k)}$ can be inductively obtained as the product of $a_1^k$ with a known functions of $\alpha_1$ and the parameters of the equation \eqref{eq_rev}.

To compute the solution for $z<0$ the procedure is perfectly the same and the expressions of the coefficients $b_k$ are formally the same of $a_k$ once substituted $\alpha_1$ with $\alpha_4=-\alpha_1$.
For example, at $k=2$ the equation one obtains:

\begin{equation}\label{step1g}
p(2\alpha_4)b_2=\left((p+q)\alpha_4^3+(c+d)\alpha_4^2+r\alpha_4+g\right)\,b_1^2,
\end{equation}

\noindent which gives $b_2$ in terms of $b_1$, in the form given in \eqref{coeff_relB}, where:

\begin{equation}\label{phi2g}
\psi_2=\frac{(p+q)\alpha_4^3+(c+d)\alpha_4^2+r\alpha_4+g}{p(2\alpha_4)}.
\end{equation}

\noindent These expressions involve both odd and even powers of $\alpha_4$, therefore there is no
longer any symmetry linking the functions $\varphi_k$ and $\psi_k$ (in the reversible case they were perfectly coincident $\varphi_k=\psi_k, k>1$, because they involved only even powers of $\alpha_1$ and $\alpha_4$).

\end{document}